\input{style/arxiv-ba.cfg}
\documentclass[ba,linksfromyear,preprint]{imsart}
\makeatletter
   \@ifpackageloaded{natbib}{}{\usepackage{natbib}}
\makeatother

\pubyear{2015}
\volume{10}
\issue{2}
\firstpage{501}
\lastpage{504}
\doi{10.1214/15-BA942A}

\begin{document}

\begin{frontmatter}
\title{Comment on Article by Dawid and Musio\thanksref{T1}}
\runtitle{Comment on Article by Dawid and Musio}
%
%
\relateddois{T1}{Main article DOI: \relateddoi[ms=BA942]{Related item:}{10.1214/15-BA942}.}

\begin{aug}
\author[a]{\fnms{Matthias} \snm{Katzfuss}\ead[label=e1]{katzfuss@tamu.edu}}
\and
\author[b]{\fnms{Anirban} \snm{Bhattacharya}\ead[label=e2]{anirbanb@stat.tamu.edu}}
\runauthor{M. Katzfuss and A. Bhattacharya}
\address[a]{Department of Statistics, Texas A\&M University,
\printead{e1}}
\address[b]{Department of Statistics, Texas A\&M University,
\printead{e2}}

\end{aug}




\end{frontmatter}

The authors consider the interesting and important issue of Bayesian
inference based on objective functions other than the likelihood. They
focus on model selection in the low-dimensional setting using
prequential local proper scoring rules.

\section{General non-likelihood-based inference}\label{sec:litreview}

There is a large and disparate literature on inference based on
objective functions other than the likelihood. We will briefly mention
some examples here, but we believe that a more thorough review and
comparison would be a worthy endeavor.

Numerous objective functions have been proposed to replace the \xch{(log-)likelihood}{(log-) likelihood}
in pursuit of various inference goals. Proper scoring rules
are a natural choice for serving as such objective functions, due to
their property of being minimized (in expectation) under the true
model. Depending on the goal of the analysis, certain well-known proper
scoring rules can achieve robustness (e.g., continuous ranked
probability score, or CRPS), have simple closed-form expressions (e.g.,
Dawid--Sebastiani score), or do not require densities (e.g., CRPS) or
normalizing constants (e.g., Hyv\"arinen score, as in the present
paper). See \citet{Gneiting2014} for a recent review of these and other
scoring rules.

In a frequentist context, examples of approaches falling into this
category of scoring-rule-based inference are minimum contrast
estimation \citep[e.g.,][]{Pfanzagl1969, Birge1993}, composite
likelihood \citep[e.g.,][]{Lindsay1988a}, and M-estimation
\citep[e.g.,][]{Huber2009}. Some further review is given in
\citet{Dawid2014}.

There have also been related approaches in the Bayesian framework.
\citet{Shaby2014} provides a nice review of Bayesian inference using
general objective functions and, based on results of
\citet{Chernozhukov2003}, he proposes an ``open-faced sandwich
adjustment'' to obtain pseudo-posteriors with properly calibrated
frequentist properties. Further, the ``Gibbs posterior''
\citep{Jiang2008,Li2013} has received considerable interest, where the
negative log-likelihood is replaced by some ``empirical risk'' $R_n$
(usually targeting the specific parameter to be estimated) to construct
a pseudo-posterior of the form
\begin{equation}
\label{eq:gibbsposterior}
Q(\theta) \propto \exp\{- \lambda R_n(\theta) \} \pi(\theta),
\end{equation}
where $\lambda$ is a positive scaling constant  (often called
``temperature''). Sampling from the pseudo-posterior $Q$ can be
performed via standard MCMC algorithms.

\section{Objective Bayesian model selection}

In objective Bayesian model selection, a discrete prior is assumed on a
(finite) class of models, and given a particular model, objective
improper priors are placed on the model parameters. While improper
priors are commonly used for analysis of a single model, one faces
difficulties in comparing models via Bayes factors, since the marginal
likelihoods of the competing models are only specified up to arbitrary
constants. A~number of remedies have been proposed in the literature to
deal with this issue, such as fractional Bayes factors
\citep{OHagan1995} and intrinsic Bayes factors \citep{Berger1996}.

In the present paper, the authors take a different approach, which
relies on replacing the \xch{(log-)marginal}{(log-) marginal} likelihood by a local proper
scoring rule. The Hyv\"arinen score is recommended as a default. From
the expression of the Hyv\"arinen score in the authors' equation (16),
it can be seen that the arbitrary constant disappears. The authors look
at examples where the Hyv\"arinen scores are analytically tractable and
provide asymptotic orders for the difference in Hyv\"arinen scores
assuming the respective models to be true.

Some clarification regarding practical implementation of the model
selection procedure presented here would be helpful. When can we be
sure that one model is truly better than another --- or in other words,
can anything be said about posterior model probabilities (also see
Section \ref{sec:scaling} below)? Can the the necessary quantities be
computed for models beyond the simple Gaussian examples considered in
the paper?

\section{Scaling issues}\label{sec:scaling}

As indicated in \eqref{eq:gibbsposterior} above, the literature on
Gibbs posteriors typically includes a multiplicative scaling constant
$\lambda$ on the objective function. The choice of $\lambda$ is
considered a critical issue, as it has a direct effect on the \xch{(pseudo-)posterior}{(pseudo) posterior}
uncertainty. \citet{Shaby2014} does not consider a
multiplicative scaling of the objective function, but his
open-face-sandwich correction automatically adjusts for such scaling,
and his approach is thus invariant to scaling. Without such a
correction, the scaling issue also arises when the objective function
is specified to be a proper scoring rule, including the Hyv\"arinen
score. As implicitly acknowledged by the authors
\xch{in their Footnote~2,}{in their second footnote on p.~6,} the scaling of a proper scoring rule is arbitrary, in
that any proper scoring rule is still proper when multiplied by a
constant.

In the context of model selection between models $M_1$ and $M_2$ with
scores $S_{M_1}$ and $S_{M_2}$, respectively, the scaling constant can
arbitrarily inflate or deflate the pseudo Bayes factor,
\[
\textnormal{PBF} = \frac{\exp(\lambda S_{M_1})}{\exp(\lambda S_{M_2})}
\]
and thus the amount of evidence in favor of $M_1$ over $M_2$
\citep[cf.][]{Kass1995}. This also makes it challenging to compute
pseudo posterior model probabilities, such as
\begin{equation}
\label{eq:ppp}
\widetilde{P}(M_1 | \mathbf{x}) = \frac{\exp(\lambda S_{M_1})}{\exp(\lambda S_{M_1}) + \exp(\lambda S_{M_2})}.
\end{equation}
If $S_{M_1}$ is larger than $S_{M_2}$, \eqref{eq:ppp} can be
arbitrarily close to $0.5$ or $1$ by choosing $\lambda$ to be very
small or very large, respectively.

In light of these scaling issues, how should model selection be
calibrated and interpreted? Moreover, is it possible to handle more
than two competing models or even high-dimensional settings, where the
number of competing models may grow exponentially with the sample size?
In the high-dimensional linear regression context,
\citet{johnson2012bayesian} showed that a number of commonly used
procedures (including fractional and intrinsic Bayes factors) assign
vanishingly small posterior probabilities to the true model with
increasing sample size. The scaling issue may assume an even more
important role in such cases.

\end{document}